\documentstyle[amssymb,12pt]{amsart}

\ifx\mathrm\undefined\let\mathrm\bf\fi
\ifx\mathbf\undefined\let\mathbf\bf\fi

 \def
\Sum{\sum\limits}

\let\dl\delta \let\Dl\Delta
\let\epe\epsilon \let\eps\varepsilon \let\epsilon\eps

\let\la\lambda \let\La\Lambda
 
 \let\phi\varphi

\newcommand{\G}{{{\frak G}\,}}
\newcommand{\F}{{{\frak F}\,}}
\newcommand{\Zb}{{{\mathcal Z}}}

\newcommand{\half}{\frac12}
\newcommand{\Z}{{\Bbb Z}}  
\newcommand{\N}{{\Bbb N}}

\newcommand{\C}{{\Bbb C}}
   
\newcommand{\Ref}[1]{{$($\ref{#1}$)$}}
\newcommand{\bean}{\begin{eqnarray}}
\newcommand{\eean}{\end{eqnarray}}
\newcommand{\be}{\begin{displaymath}}
\newcommand{\ee}{\end{displaymath}}
\newcommand{\bea}{\begin{eqnarray*}}   
\newcommand{\eea}{\end{eqnarray*}}
\newcommand{\g}{{{\frak g}\,}}

\newcommand{\Imag}{{\operatorname{Im}}\,}
\newcommand{\Real}{{\operatorname{Re}}\,}

\newcommand{\T}{\otimes}

\newcommand{\vs}{\vspace{1.5\baselineskip}}

\newenvironment{proof}{\noindent{\it Proof\/}:\rm}{$\;\Box$
\par\vs}

\newtheorem
{thm}{Theorem}

\newtheorem
{lemma}[thm]{Lemma}

\newtheorem
{corollary}[thm]{Corollary}

\newcommand{\th}{\theta}
\newcommand{\End}{{\operatorname{End}}}

\begin{document}

\title[On algebraic equations satisfied by hypergeometric solutions...]{On algebraic equations satisfied by hypergeometric solutions of the qKZ equation.} 
\author{E. Mukhin and A. Varchenko}
\maketitle
\vskip-.5\baselineskip
\centerline{{\it Department of Mathematics,
University of North Carolina at Chapel Hill,}}
\centerline{\it Chapel Hill, NC 27599-3250, USA}
\centerline{{\it E-mail addresses:} {\rm mukhin@@math.unc.edu,
av@@math.unc.edu}}
\bigskip
\medskip
\centerline{October, 1997}
\bigskip
\medskip

\begin{abstract}
We consider the $sl(2)$ quantized Knizhnik-Zamolodchikov equation (qKZ), defined in terms of rational R-matrices. The properties of the equation change when the step of the equation takes a resonance value. In this case the discrete connection defined by the qKZ equation has a invariant subbundle which we call the subbundle of quantized conformal blocks. Solutions of the qKZ equation were constructed in \cite{TV1},\cite{MV1} in terms of multidimensional hypergeometric integrals. In this paper we show that for a resonance step all  hypergeometric solutions take values in the subbundle of quantized conformal blocks, moreover the values span the subbundle of quantized conformal blocks under certain conditions. We describe the space of hypergeometric solutions in terms of the quantum group $U_qsl(2)$.
\end{abstract}

\section{Introduction}

Conformal field theory associates a finite dimensional vector space, called the space of conformal blocks, to each Riemann surface with marked points and certain additional data (local coordinates, representations). The vector spaces corresponding to different complex structures or different positions of the points are locally (projectively) identified by a projectively flat connection.

A Wess-Zumino-Witten model is labeled by a simple Lie algebra $\g$ and a positive integer $c$ called level. The space of conformal blocks is defined in terms of the representation theory of the affine Kac-Moody algebra $\widehat{L\g}$, which is a central extension of the loop algebra $L\g=\g\T\C((t))$, see \cite{K}. For each irreducible highest weight $\g$-module $L$, we have a canonically defined corresponding irreducible highest weight $\widehat{L\g}$ module of level $c$, denoted $\widehat{L}$, see \cite{K}. Then, the space of conformal blocks associated to a Riemann surface $\Sigma$, $n$ distinct points with a choice of local holomorphic parameters around them and $n$ irreducible finite dimensional $\g$ modules $L_1,\dots, L_n$ (such that $\widehat{L}_i$ are integrable) is the space $H^0(\g(\Sigma),(\widehat{L}_1\T\ldots\T \widehat{L}_n)^*)$ of linear forms on $\widehat{L}_1\T\ldots\T \widehat{L}_n$ invariant under the action of the Lie algebra $\g(\Sigma)$ of meromorphic $\g$-valued functions on $\Sigma$, holomorphic outside of the marked points. The action of $\g(\Sigma)$ is defined through the Laurent expansion at the poles. Varying the data makes the spaces of conformal blocks into a holomorphic vector bundle with a projectively flat connection given by the Sugawara-Segal construction.

An explicit description is known on the Riemann sphere ${\Bbb P}^1$. Namely, the space of conformal blocks on ${\Bbb P}^1$ with $n+1$ distinct points $z_1,\dots,z_n\in\C\subset{\Bbb P}^1$, $z_{n+1}=\infty$, associated to $n+1$ irreducible $\g$ modules $L_1,\dots,L_{n+1}$ with highest weights $\la_1,\dots,\la_{n+1}$, is identified with a subspace of $(L_1\T\ldots\T L_n)^{sing}_\la$, where  $(L_1\T\ldots\T L_n)^{sing}_\la\subset L_1\T\ldots\T L_n$ is the weight subspace of singular vectors of total weight $\la$, $\la=-w\la_{n+1}$ and $w$ is the longest element of the Weyl group. More precisely the space of conformal blocks can be described as follows. Let $\theta$ be the highest root and let the scalar product be normalized by $(\theta,\theta)=2$.

If the resonance condition, 
\bean\label{res'}
(\theta,\la)-c+k-1=0,
\eean
holds for some $k\in\N$, then the space of conformal blocks $N_{\la_{n+1}}(z)$ is identified with the subspace
%$k\le (\th,\la_{n+1})$,
\be
N_{\la_{n+1}}(z)=\{m\in (L_1\T\ldots\T L_n)^{sing}_\la\,|\, (E(z))^km=0\},
\ee
otherwise the space of conformal blocks is identified with the weight space of singular vectors,
$N_{\la_{n+1}}(z)=(L_1\T\ldots\T L_n)^{sing}_\la$.
Here
\be
E(z)=\sum_{i=0}^n z_ie_\th^{(j)},
\ee
and $e_\th^{(j)}$ denotes the element $e_\th\in \g$ acting on the $j$-th factor, see \cite{KT}, \cite{FSV1}.

The subspaces $N_{\la_{n+1}}(z)$ for different sets of $z_1,\dots,z_n\in\C$ form a subbundle of the trivial vector bundle over the space $\C^n$ with fiber $L_1\T\ldots\T L_n$. There is a flat connection on the bundle $\C^n\times (L_1\T\ldots\T L_n)$ preserving the subbundle of conformal blocks. Its horizontal sections $\Psi(z)$ obey the Knizhnik-Zamolodchikov equation, 
\bean\label{KZ}
\partial_{z_i}\Psi(z)=\frac{1}{\kappa}\sum_{j\neq i}\frac{\Omega_{ij}}{z_i-z_j}\Psi(z),\qquad i=1,\dots,n,
\eean
$\kappa=c+h^\vee$, where $\Omega_{ij}$ is the Casimir operator acting on the $i$-th and $j$-th factors and $h^\vee$ is the dual Coxeter number of $\g$.

In \cite{SV} the KZ equation was solved in terms of hypergeometric integrals, cf. \cite{R}.
It was shown in \cite{V} that for generic values of $\kappa$, the space of hypergeometric solutions can be naturally identified with the space of singular vectors in the tensor product $L_1^q\T\ldots\T L_n^q$ of the corresponding representations of the quantum group $U_q\g$ with $q=e^{-\pi i/\kappa}$.    

In [FSV1-3] it was shown that under resonance condition \Ref{res'}, all hypergeometric solutions of the KZ equation with values in $(L_1\T\ldots\T L_n)_\la^{sing}$ automatically take values in the space of conformal blocks, i.e. their values lie in the kernel of the operator $(E(z))^k$. Moreover, in \cite{V} it was shown that for $\g=sl(2)$, the hypergeometric solutions span the space of conformal blocks and the space of solutions is naturally identified with the co-image of the operator $(f_q)^k\in U_qsl(2)$,
\be
(L_1^q\T\ldots\T L_n^q)_{\la}^{sing}/(f_q)^k((L_1^q\T\ldots\T L_n^q)_{\la+k\th}^{sing}).
\ee

Let $\g=sl(2)$. The quantized Knizhnik-Zamolodchikov (qKZ) equation is a holonomic system of difference equations for a function $\Psi(z)$ with values in $L_1\T\ldots\T L_n$,
\be
\Psi(z_1,\dots,z_i+p,\dots,z_n)=R_{m,m-1}(z_m-z_{m-1}+p)\ldots R_{m,1}(z_m-z_1+p)\times
\ee
\be
\times
R_{m,n}(z_m-z_n)\ldots R_{m,m+1}(z_m-z_{m+1})\Psi(z),
\ee
$i=1,\dots,n$, where $R_{i,j}(x)$ is the rational $R$-matrix $R_{L_iL_j}(x)$ acting in $i$-th and $j$-th factors and $p\in\C$ is a parameter, see \cite{FR}.

The qKZ equation defines the discrete Knizhnik-Zamolodchikov (qKZ) connection on the trivial vector bundle over $\C^n$ with fiber $L_1\T\ldots\T L_n$.
 
Consider the quasiclassical limit of the qKZ. Namely, set $y_i=z_i/h$, for some $h\in\C$ and let $h\to 0$. In this limit the qKZ equation turns into a system of differential equations
\bean\label{KZ'}
p\partial_{y_i}\tilde{\Psi}(y)=-\sum_{j\neq i}\frac{\tilde{\Omega}_{ij}}{y_i-y_j}\tilde{\Psi}(y),\qquad i=1,\dots,n,
\eean
where $\tilde{\Omega}_{ij}=\Omega_{ij}-2\la_i\la_j$ see Section 7 in \cite{TV1} and Section 12.5 in \cite{CP}.

Notice that if $p=-\kappa$, then $\Psi(z)$ is a solution of the KZ equation \Ref{KZ} if and only if the function
\be
\tilde{\Psi}(y)=\prod_{i<j}(y_i-y_j)^{-2\la_i\la_j/\kappa}\,\Psi(y)
\ee
is a solution of the equation \Ref{KZ'}.

In \cite{TV1}, \cite{MV1} the qKZ equation was solved in terms of hypergeometric integrals. It was shown that the space of hypergeometric solutions for generic values of $p$ can be naturally identified with the space of singular vectors in the tensor product $L_1^q\T\ldots\T L_n^q$ of the corresponding representations of the quantum group $U_qsl(2)$ with $q=e^{\pi i/p}$.    

In \cite{MV2} a quantization of the space of conformal blocks is suggested. Namely, under the resonance condition, 
\bean\label{res}
2\la+p+N+k-1=0,
\eean
where $k\in\N$, we define the space of quantized conformal blocks, $C_{\la_{n+1}}(z)$, as
%$k\le (\th,\la_{n+1})$
\be
C_{\la_{n+1}}(z)=\{m\in (L_1\T\ldots\T L_n)^{sing}_\la\,|\, (e(z))^km=0\},
\ee
where
\be
e(z)m=\sum_{j=1}^n \left(z_j+h^{(j)}+\sum_{s=j+1}^n2h^{(s)}\right)e^{(j)}m,
\ee
and $h^{(s)},e^{(s)}$ denote the elements $h,e\in sl(2)$ acting on the $s$-th factor of $L_1\T\ldots\T L_n$. 

The operator $e(z)$ can be described in terms of the action of the Yangian $Y(gl(2))$ in the tensor product of evaluation modules $L_1(z_1)\T\ldots\T L_n(z_n)$, $e(z)=T_{21}^{(2)}-T_{22}^{(1)}T_{21}^{(1)}$, see \Ref{e(z)-Yangian}. In the quasiclassical limit the resonance condition \Ref{res} coincides with the resonance condition \Ref{res'} and the operator $e(z)$ tends to the operator $E(z)$.

It was shown in \cite{MV2} that the spaces of quantized conformal blocks form a subbundle invariant with respect to the quantized KZ connection.

In this paper we show that under resonance condition \Ref{res}, all hypergeometric solutions of the $sl(2)$ qKZ equation with values in $(L_1\T\ldots\T L_n)^{sing}_\la$ automatically take values in the space of quantized conformal blocks. We identify the space of the hypergeometric solutions with the co-image of the operator $(f_q)^k\in U_qsl(2)$,
\be
(L_1^q\T\ldots\T L_n^q)_{\la}^{sing}/(f_\th^q)^k((L_1^q\T\ldots\T L_n^q)_{\la+k\th}^{sing}).
\ee
We prove that under certain conditions the hypergeometric solutions span the space of quantized conformal blocks.

The paper is organized as follows. In Section~\ref{general} we recall some general facts about $sl(2)$, $U_qsl(2)$ and $Y(gl(2))$. In Section~\ref{hypergeometric} we recall the construction of the hypergeometric solutions of the qKZ equation. In Section~\ref{s(z)} we state and prove the main results of this paper.

\section{General definitions and notations}\label{general}
\subsection{The Lie algebra $sl(2)$}

Let $e$, $f$, $h$ be generators of the Lie algebra $sl(2)$ such that 
\be
[h,e]=e, \qquad [h,f]=-f, \qquad [e,f]=2h.
\ee
For an $sl(2)$ module $M$, 
let $M^*$ be its restricted dual with an $sl(2)$ module structure defined by
\be
\langle e\varphi,x\rangle=\langle \varphi,fx\rangle,  \qquad
\langle f\varphi,x\rangle=\langle \varphi,ex\rangle,  \qquad
\langle h\varphi,x\rangle=\langle \varphi,hx\rangle
\ee
for all $x\in M$, $\varphi \in M^*$. The module $M^*$ is called the \emph{dual} 
module. 

For $\la \in \C$, 
denote $V_\la$ the $sl(2)$  Verma module with  highest weight $\la$,
$V_\la=\bigoplus\limits_{i=0}^\infty \C f^iv$, where $v$ is a highest weight
vector. Denote $L_\la$ the irreducible module with highest weight $\la$.

Let $\La^+=\{0,\frac{1}{2},1,\frac{3}{2},2,...\}$ be the set of dominant 
weights. If $\la\in\La^+$, then $L_\la$ is  a $(2\la+1)$-dimensional module and
\be
L_\la\,\simeq\, V_\la/S_\la,
\ee
where $S_\la =\bigoplus\limits_{i=2\la+1}^\infty \C f^iv\subset
V_\la$ is the 
maximal proper submodule.  The
vectors $f^iv$, $i=0,\dots,2\la$, generate a
basis in $L_\la$.

For $\la\not\in\La^+$, $L_\la =V_\la$.
It is convenient to introduce $S_\la$ to be the zero submodule of $V_\la$, then
$L_\la\simeq V_\la/S_\la$ for all $\la\in\C$.

For an $sl(2)$ module $M$ with  highest weight $\la$, denote by $(M)_l$ the 
subspace of weight $\la-l$, by $(M)^{sing}$ the kernel of the operator $e$, and by $(M)_l^{sing}$ the subspace $(M)_l\bigcap (M)^{sing}$.

\subsection{The algebra $U_qsl(2)$.} 

Let q be a complex number different from $\pm 1$. 
Let $e_q$, $f_q$, $q^h$, $q^{-h}$ be generators of $U_qsl(2)$ such that
\be
q^hq^{-h}=q^{-h}q^h=1, \qquad [e_q,f_q]=\frac{q^{2h}-q^{-2h}}{q-q^{-1}},
\ee
\be
q^he_q=qe_qq^h, \qquad  q^hf_q=q^{-1}f_qq^h.
\ee

A comultiplication $\Delta : U_qsl(2)\to U_qsl(2)\T U_qsl(2)$ is given by
\be
\Delta(q^h)=q^h\T q^h,\qquad \Delta(q^{-h})=q^{-h}\T q^{-h},
\ee
\be
\Delta(e_q)=e_q\T q^h + q^{-h}\T e_q,\qquad 
\Delta(f_q)=f_q\T q^h + q^{-h}\T f_q.
\ee
The comultiplication defines a module structure on tensor products of 
$U_qsl(2)$ modules.

For $\la \in \C$, denote $V_\la^q$ the $U_qsl(2)$ Verma module with highest weight 
$q^\la$, $V^q_\la=\bigoplus\limits_{i=0}^\infty \C f_q^iv^q$, where $v^q$ is a highest weight
vector. 

For $\la\in \La^+$, $S^q_\la=\bigoplus\limits_{i=2\la+1}^\infty \C f_q^iv^q$ is a submodule in $V^q_\la$.
Denote $L^q_\la$ the quotient module $V^q_\la/S^q_\la$. The module
$L^q_\la$ is the $(2\la+1)$-dimensional  
highest weight  module with highest weight $q^\la$. The 
vectors $f_q^iv^q$, $i=0,\dots,2\la$, generate a 
basis in $L^q_\la$.

For $\la\not\in\La^+$, let $L^q_\la=V^q_\la$.
It is convenient to introduce $S^q_\la$ to be the zero submodule of $V^q_\la$, then
$L^q_\la\simeq V^q_\la/S^q_\la$ for all $\la\in\C$.

For an $U_qsl(2)$ module $M^q$ with highest weight $q^\la$, denote by 
$(M^q)_l$ the
subspace of weight $q^{\la-l}$, by $(M^q)^{sing}$ the kernel of the 
operator $e_q$, and
by $(M^q)_l^{sing}$ the subspace $(M^q)_l\bigcap (M^q)^{sing}$.

\subsection{The Hopf algebra $Y(gl(2))$.}
The Yangian $Y(gl(2))$ is an associative algebra with an infinite set of generators $T_{i,j}^{(s)}$,  $i,j=1,2$, $s=0,1,\dots$, subject to the following relations:
\be
[T_{ij}^{(r)},T_{kl}^{(s+1)}]-[T_{ij}^{(r+1)},T_{kl}^{(s)}]=T_{kj}^{(r)}T_{il}^{(s)}-
T_{kj}^{(s)}T_{il}^{(r)},\qquad T_{ij}^{(0)}=\dl_{ij},
\ee
$i,j,k,l=1,2$; $r,s,=1,2,\dots$ .

The comultiplication $\Dl\, :\,Y(gl(2))\to Y(gl(2))\T Y(gl(2))$ is given by 
\be
\Dl\,:\,T_{ij}^{(s)}\mapsto \sum_{k=1}^2\sum_{r=0}^s\, T_{ik}^{(r)}\T T_{kj}^{(s-r)}.
\ee 

For each $x\in\C$, there is an automorpfism $\rho_x\,:\, Y(gl(2))\to Y(gl(2))$ given by
\be
\rho_x\,:\,T_{ij}^{(s)}\mapsto\sum_{r=1}^s\,{s-1 \choose r-1}x^{s-r}T_{ij}^{(r)}.
\ee

There is also an \emph{evaluation morphism} $\epe$ to the universal enveloping algebra of $sl(2)$,
$\epe\,:\, Y(gl(2))\to U(sl(2))$, given by
\be
\epe\,:\,T_{11}^{(s)}\mapsto \dl_{1s}h,\qquad \epe\,:\,T_{12}^{(s)}\mapsto \dl_{1s}f,
\ee
\be
\epe\,:\,T_{21}^{(s)}\mapsto \dl_{1s}e,\qquad \epe\,:\,T_{22}^{(s)}\mapsto -\dl_{1s}h,
\ee
for $s=1,2,\dots$ .

Introduce the generating series $T_{ij}(u)=\Sum_{s=0}^\infty T_{ij}^{(s)}u^{-s}$. In terms of these series the relations in the Yangian take the form
\be
R(x-y)T_{(1)}(x)T_{(2)}(y)=T_{(2)}(y)T_{(1)}(x)R(x-y),
\ee
where $R(x)=(x\operatorname{Id}+P)\in\End(\C^2\T\C^2)$, $P\in\End(\C^2\T\C^2)$ is the operator of permutation of the two factors, 
$T_{(1)}(x)=1\T T(x)$, $T_{(2)}(x)=T(x)\T 1$.

In terms of the generating series the comultiplication $\Dl$, the automorphisms $\rho_x$ and the evaluation morphism $\epe$ take the form
\be
\Dl\,:\, T_{ij}(u)\mapsto \sum_{k=1}^2 T_{ik}(u)\T T_{kj}(u),
\ee
\be
\rho_x\,:\,T_{ij}(u)\mapsto T_{ij}(u-x),
\ee
\be
\epe\,:\,T_{11}(u)\mapsto \frac{h}{u},\qquad \epe\,:\,T_{12}(u)\mapsto \frac{f}{u},
\ee
\be
\epe\,:\,T_{21}(u)\mapsto \frac{e}{u},\qquad \epe\,:\,T_{22}(u)\mapsto -\frac{h}{u},
\ee
$i,j=1,2$. For more detail on the Yangian see \cite{CP},\cite{KR}.

For any $sl(2)$ module $M$ and $x\in\C$, let $M(x)$ be the $Y(gl(2))$ module obtained from $M$ via the homomorphism $\epe\circ\rho(x)$. The module $M(x)$ is called the \emph{evaluation module}. The action of $Y(gl(2))$ in the evaluation module $M(x)$ is given by
\be
T_{11}^{(s)}m=x^{s-1}h\,m,\qquad T_{12}^{(s)}m=x^{s-1}f\,m,
\ee 
\be
T_{21}^{(s)}m=x^{s-1}e\,m,\qquad T_{22}^{(s)}m=-x^{s-1}h\,m,
\ee 
for all $m\in M$, $s=1,2,\dots$ .

For $\la_1,\la_2\in\C$ and generic complex numbers $x,y$, the $Y(gl(2))$ modules $L_{\la_1}(x)\T L_{\la_2}(y)$ and $L_{\la_2}(y)\T L_{\la_1}(x)$ are irreducible and isomorphic. There is a unique intertwiner of the form $PR_{L_{\la_1}L_{\la_2}}(x-y)$ mapping $v_1\T v_2$ to $v_2\T v_1$, where $P$ is the operator of permutation of the two factors and $v_i$ are highest weight vectors generating $L_{\la_i}$, $i=1,2$. The operator $R_{L_{\la_1}L_{\la_2}}(x)\in\End(L_{\la_1}\T L_{\la_2})$ is called the \emph{rational $R$-matrix}, see  \cite{CP}, \cite{D}.

The vector spaces $V_{\la_1}\T V_{\la_2}$ for different values of $\la_1,\la_2\in\C$ are identified by distinguished bases $\{f^{l_1}v_1\T f^{l_2}v_2\,|\,l_1,l_2\in\Z_{\ge 0}\}$.

\begin{thm}\label{R-matrix}(Theorem 1 in \cite{MV1}.) 

1. The rational $R$-matrix $R_{V_{\la_1}V_{\la_2}}(x)\in\End(V\T V)$ is a meromorphic function of $x,\la_1,\la_2$. Moreover, for any $\la_1,\la_2\in\C$, the $R$-matrix $R_{V_{\la_1}V_{\la_2}}(x)$ is a well defined meromorphic function of $x$.

2. For generic $x\in\C$, the rational $R$-matrix $R_{V_{\la_1}V_{\la_2}}(x)$ preserves the submodule $V_{\la_1}\T S_{\la_2}+S_{\la_1}\T V_{\la_2}\subset V_{\la_1}\T V_{\la_2}$.

3. Let $V_{\la_1}\T V_{\la_2}\to L_{\la_1}\T L_{\la_2}$ be the canonical factorization map. Then, for generic $x$,
the rational $R$-matrix $R_{V_{\la_1}V_{\la_2}}(x)$ can be factorized to an operator $R(x)\in\End(L_{\la_1}\T L_{\la_2})$ and, moreover, $R(x)=R_{L_{\la_1}L_{\la_2}}(x)$.
\end{thm}

Let either $M_i=V_{\la_i}$, $i=1,2$, be Verma $sl(2)$ modules or $M_i=L_{\la_i}$, $i=1,2$, be irreducible $sl(2)$ modules. For all $x\in\C$, the rational $R$-matrix $R_{M_1M_2}(x)$ commutes with the action of $sl(2)$ in $M_1\T M_2$ and, in particular, preserves the weight decomposition.

\subsection{The qKZ connection.}\label{qKZ}

The \emph {rational quantized Knizhnik-Zamolodchikov equation (qKZ)} associated 
to $sl(2)$ is the following holonomic system of linear
difference equations for a
function $\Psi(z_1,\dots,z_n)$ with values in a tensor product 
$M_1\T\ldots\T M_n$ of $sl(2)$ modules: 
\be
\Psi(z_1,\dots,z_m+p,\dots,z_n)=K_m(z)\Psi(z_1,\dots,z_n),\qquad m=1,\dots,n,
\ee
\bea
\lefteqn{K_m(z)=R_{M_mM_{m-1}}(z_m-z_{m-1}+p)\dots R_{M_mM_1}(z_m-z_1+p)\times}
\\&&
\times R_{M_mM_n}(z_m-z_n)\dots R_{M_mM_{m+1}}(z_m-z_{m+1}),
\eea
where $p$ is a complex parameter, $R_{M_iM_j}(x)\in \End(M_i\T M_j)$ is the rational $R$-matrix acting in the $i$-th and $j$-th factors, see \cite{FR}. The linear operators $K_i(z)$ are called the \emph{qKZ operators}.

The qKZ operators commute with the $sl(2)$ action in the tensor product $M_1\T\ldots\T M_n$ and, in particular, preserve the subspaces $(M_1\T\ldots\T M_n)_l$ and $(M_1\T\ldots\T M_n)_l^{sing}$ for all $l\in\Z_{\ge 0}$. In order to construct all solutions of the qKZ equation, it is enough to solve the equation with values in singular weight spaces $(M_1\T\ldots\T M_n)_l^{sing}$.

The qKZ operators define a discrete flat connection on the trivial vector bundle over $\C^n$ with fiber $M_1\T\ldots\T M_n$. This connection is called the \emph{quantized Knizhnik-Zamolodchikov connection}.

\section{Hypergeometric solutions of the qKZ equation.}\label{hypergeometric}

\subsection{The phase function}

Let $z=(z_1,\dots,z_n)\in\C^n,\;\la=(\la_1,\dots,\la_n)\in\C^n ,\; 
t=(t_1,\dots,t_l)\in\C^l$.
The \emph{phase function} is defined by the following formula:
\be
\Phi(t,z,\la)=
\prod_{i=1}^n\prod_{j=1}^l
\frac{\Gamma((t_j-z_i+\la_i)/p)}{\Gamma((t_j-z_i-\la_i)/p)}
\prod_{1\le i< j\le l}\frac{\Gamma((t_i-t_j-1)/p)}{\Gamma((t_i-t_j+1)/p)}.
\ee

\subsection {Actions of the symmetric group}\label{actions}

Let $f=f(t_1,\dots,t_l)$ be a function. For a permutation $\sigma\in{\Bbb S}^l$, define the functions
$[f]_{\sigma}^{rat}$ and $[f]_{\sigma}^{trig}$ via the action of the simple transpositions 
$(i,i+1)\in {\Bbb S}^l$, $i=1,\dots,\l-1$, given by      

\be
[f]_{(i,i+1)}^{rat}(t_1,\dots, t_l)=
f(t_1,\dots,t_{i+1},t_i,\dots,t_l)\frac{t_i-t_{i+1}-1}{t_i-t_{i+1}+1},
\ee
\be
[f]_{(i,i+1)}^{trig}(t_1,\dots, t_l)=
f(t_1,\dots,t_{i+1},t_i,\dots,t_l)
\frac{
\sin(\pi(t_i-t_{i+1}-1)/p)}{\sin(\pi(t_i-t_{i+1}+1)/p)}.
\ee

If for all $\sigma\in{\Bbb S}^l$, $[f]_{\sigma}^{rat}=f$, then we say that the function is
\emph{symmetric with respect to the rational action}.
If for all $\sigma\in{\Bbb S}^l$, $[f]_{\sigma}^{trig}=f$, then we say that the function is
\emph{symmetric with respect to the trigonometric action}.
 
This definition implies the following important Remark.

{\bf Remark.} If $w(t_1,\dots, t_l)$ is ${\Bbb S}^l$ symmetric 
with respect to the 
rational action and $W(t_1,\dots, t_l)$ is ${\Bbb S}^l$ symmetric with respect to the 
trigonometric action,
then $\Phi wW$ is a symmetric function  of $t_1,\dots, t_l$ (in the 
usual sense).

\subsection {Rational weight functions}

Fix natural numbers $n,l$.

Set $\Zb^n_l=\{\bar{l}=(l_1,\dots,\l_n)\in\Z^n_{\ge 0}\,|\, 
\sum\limits_{i=1}^nl_i=l\}$. For $\bar{l}\in\Zb^n_l$ 
and $m=0,1,\dots,n$, set $l^m=\sum\limits_{i=1}^ml_i$.

For $\bar{l}\in\Zb_l^n$, define the \emph{rational weight function} 
$w_{\bar{l}}$ by $w_{\bar{l}}=\Sum_{\sigma\in {\Bbb S}^l}[\eta_{\bar{l}}]^{rat}_\sigma$, where
\bean\label{eta}
\eta_{\bar{l}}(t,z,\la)=  
\prod_{m=1}^n\frac{1}{l_m!}\prod_{j=l^{m-1}+1}^{l^m}
\left( \frac{1}{t_j-z_m-\la_m}\prod_{k=1}^m\frac{t_j-z_k+\la_k}{t_j-z_k-\la_k}\right).
\eean

{\bf Example.} Let $l=1$. Then the rational weight functions $w_{(0,\dots,1_i,\dots,0)}=w_i$ have the form
\bean\label{one point}
w_i(t,z,\la)=\frac{1}{t-z_i-\la_i}\prod_{m=1}^{i-1}\frac{t-
z_m+\la_m}{t-z_m-\la_m},
\eean
$i=1,\dots,n$.

\medskip

For fixed $z,\la\in\C^n$, the space spanned over $\C$
by  all rational weight functions 
$w_{\bar{l}}(t,z,\la)$, $\bar{l}\in\Zb_l^n$, is called 
the \emph{hypergeometric rational space specialized at $z,\la$} and 
is denoted $\F(z,\la)=\F^n_l(z,\la)$.
This space is a space of functions of variable $t$.

\subsection{Trigonometric weight functions}\label{trig}
Fix natural numbers $n, l$.

For $\bar{l}\in\Zb^n_l$, 
define the \emph{trigonometric weight function} $W_{\bar{l}}$ by
\bea
\lefteqn{W_{\bar{l}}(t,z,\la)=}
\\
&&
\sum_{\sigma \in {\Bbb S}^l}
\left[\prod_{m=1}^n\prod_{s=1}^{l_m}
\frac{\sin(\pi/p)}{\sin(\pi s/p)}
\prod_{j=l^{m-1}+1}^{l^m}\!\!
\frac{\exp(\pi i(z_m-t_j)/p)}{\sin(\pi(t_j-z_m-\la_m)/p)}
\prod_{k=1}^m\frac{\sin(\pi(t_j-z_k+\la_k)/p)}{\sin(\pi(t_j-z_k-\la_k)/p)}
\right] _\sigma^{trig}.
\eea

A function $W(t,z,\la)$ is said to be a \emph{holomorphic trigonometric weight function} if
\bean\label{trig.decomposition}
W(t,z,\la)=\sum_{m_1+\ldots +m_n=l}a_{\bar{m}}(\la,e^{2\pi iz_1/p},\dots,e^{2\pi iz_n/p})
W_{\bar{m}}(t,z,\la),
\eean
where $a_{\bar{m}}(\la,u)$ are holomorphic functions of parameters $\la,u\in\C^n$.
We denote $\G$ the space of all holomorphic trigonometric weight functions. This space
is a space of functions of variables $t,z,\la$.

For fixed $\la, z\in\C^n$, the space spanned over $\C$
by  all trigonometric weight functions
$W_{\bar{l}}(t,z,\la),\,\bar{l}\in\Zb_l^n$, is called
the \emph{hypergeometric trigonometric space specialized at $z,\la$} and is denoted
$\G(z,\la)=\G^n_l(z,\la)$. This space is a space of functions of variable $t$.

\medskip 

For $\bar{l}\in\Zb^{n-1}_l$, define the \emph{singular trigonometric weight function}
$W_{\bar{l}}^{sing}$ by
\bea
\lefteqn{W_{\bar{l}}^{sing}(t,z,\la)=}
\\
&&
\sum_{\sigma \in {\Bbb S}^l}
\left[\prod_{m=1}^{n-1}\prod_{s=1}^{l_m}
\frac{\sin(\pi/p)}{\sin(\pi s/p}
\sin(\pi (z_m-\la_m-z_{m+1}-\la_{m+1}+s-1)/p)\times\right.
\\&&
\times
\prod_{j=l^{m-1}+1}^{l^m}\!\!
\frac{1}{\sin(\pi(t_j-z_m-\la_m)/p)\sin(\pi(t_j-z_{m+1}-\la_{m+1})/p)}\times
\\&&
\left.\times 
\prod_{k=1}^m\frac{\sin(\pi(t_j-z_k+\la_k)/p)}{\sin(\pi(t_j-z_k-\la_k)/p)}
\right] _\sigma^{trig}.
\eea

For fixed $z,\la\in\C^n$, the space spanned over $\C$ by all singular
trigonometric weight
functions $W_{\bar{l}}(t,z,\la),\,\bar{l}\in\Zb_l^n$,
is called the \emph{singular hypergeometric trigonometric space specialized at $z,\la$} and
is denoted $\G^{sing}(z,\la)=\G^{sing,n}_l(z,\la)$.

We have $\G^{sing}(z,\la)\subset\G(z,\la)$, see Lemma 2.29 in
\cite{TV1}.

A function $W(t,z,\la)\in\G$ is said to be a \emph{holomorphic singular trigonometric weight
function} if
$W(t,z,\la)$ is a holomorphic trigonometric weight function,
and for all $z,\la\in\C^n$, the function $W(t,z,\la)$ belongs to $\G^{sing}(z,\la)$.
We denote $\G^{sing}$ the space of all holomorphic singular trigonometric weight functions. This
space is a space of functions of variables $t,z,\la$.

For any $\bar{l}\in\Zb^{n-1}_l$, the function
$W_{\bar{l}}^{sing}(t,z,\la)$ belongs to $\G^{sing}$, see Lemma 4 in \cite{MV1}.

\medskip

Let $\la=(\la_1,\dots,\la_n)\in\C^n$, $\bar{l}=(l_1,\dots,l_n)\in\Z^n_{\ge 0}$.
An $i$-th coordinate of $\bar{l}$ is called
\emph{$\la$-admissible} if either $\la_i\not\in\La^+$
or $\la_i\in\La^+$ and $l_i\le 2\la_i$.
An index $\bar{l}$ is called \emph{$\la$-admissible} if all its coordinates are
$\la$-admissible.

For fixed $z,\la\in\C^n$, the space spanned over $\C$ by trigonometric weight functions $W_{\bar{l}}(t,z,\la)$ with $\la$-admissible indices $\bar{l}$ is called the \emph{$\la$-admissible trigonometric hypergeometric space specialized at $z,\la$} and is denoted $\G_{\rm adm}(z,\la)$.

For $\la\in\C^n$, a function $W(t,z,\mu)\in\G$ is called $\la$-\emph{admissible} if for all non-$\la$-admissible indices $\bar{m}\in\Zb^n_l$, the functions $a_{\bar{m}}(\mu,u)$ in decomposition \Ref{trig.decomposition} are equal to zero.

For fixed $z,\la\in\C^n$, the space spanned over $\C$ by all $\la$-admissible holomorphic singular trigonometric weight functions is called the \emph{$\la$-admissible singular hypergeometric trigonometric space specialized at $z,\la$} and is denoted $\G^{sing}_{\rm adm}(z,\la)$. 

\subsection{Hypergeometric integrals}\label{integrals}

Fix $p\in\C,\,\Real p<0$. Assume that the parameters $z,\la\in\C^n$ satisfy the condition $\Real (z_i+\la_i)<0$ and $\Real (z_i-\la_i)>0$ for all $i=1,\dots,n$. For a rational weight function $w=w_{\bar{l}}(t,z,\la)$, $\bar{l}\in\Zb^n_l$, and a singular trigonometric weight function $W(t,z,\la)\in\G^{sing}$, define the \emph{hypergeometric integral} $I(w,W)(z,\la)$ by the formula
\bean\label{int}
I(w,W)(z,\la)=\int\limits_{\Real t_i=0,\atop
i=1,\dots,l}\Phi(t,z,\la)w(t,z,\la)W(t,z,\la)\,d^lt,
\eean
where $d^lt=dt_1\ldots dt_l$.

The hypergeometric integral for generic $z,\la$ and an arbitrary step $p$ with negative real part is defined by analytic continuation with respect to $z,\la$ and $p$, see \cite{TV1}, \cite{MV1}.

\medskip

For a function $W(t,z,\la)\in\G^{sing}$, let $\Psi_W(z,\la)$ be the following $V_{\la_1}\T\ldots\T V_{\la_n}$-valued function
\bean\label{Psi}
\Psi_W(z,\la) \, = \, \sum_{l_1+\ldots+l_n=l} \,
I(w_{\bar{l}},W)(z,\la) \, f^{l_1}v_1\T\ldots\T f^{l_n}v_n. 
\eean

\begin{thm}\label{cite3}
(Corollaries 5.25, 5.26 in \cite{TV1}.)
Let $\la=(\la_1,\dots,\la_n)\in \C^n$ and $\la_i\not\in\La^+$.
Then for generic values of $p$ and for any function $W\in\G^{sing}$, the function
$\Psi_W(z,\la)$ is a meromorphic solution of the qKZ equation with values in $(V_{\la_1}\T\ldots\T V_{\la_n})_l^{sing}$.
\end{thm}

We always asume the following conditions on parameters $p\in\C$, $z,\la\in\C^n$:
\bean\label{step}
\Real p<0,\qquad 1\not\in p\Z,
\eean
\bean\label{weights1}
\{s\, | \, s\in\Z_{>0}, s<
2\max\{\Real\la_1,\dots,\Real\la_n\},\, s\le l\}\bigcap\{p\Z\}=\emptyset,
\eean
\bean\label{weights2}
\{2\la_m-s\,|\,s\in\Z_{\ge 0},\; s<2\Real\la_m,\, s<l\}\bigcap \{p\Z\}=\emptyset, \qquad
m=1,\dots,n,
\eean
\bean\label{resonance'}
z_k-z_m\pm(\la_k+\la_m)+s \not\in\{p\Z\},\qquad k,m=1,\dots,n,\,k\neq m,\, s=1-l,\dots,l-1.
\eean
For each $i\in\{1,\dots,n\}$ such that $\la_i\not\in\La^+$, assume
\bean\label{step3}
\{1,\dots,l\}\bigcap\{p\Z\}=\emptyset,
\eean
\bean\label{weights3}
\{2\la_i-s\,|\,s=0,1,\dots,l-1\}\bigcap \{p\Z\}=\emptyset.
\eean

Let $\la\in\C^n$. For a $\la$-admissible function $W\in\G^{sing}$, consider a function
\bean\label{Psi adm}
\Psi^{\rm adm}_W(z,\la)=\sum
I(w_{\bar{l}},W)(z,\la) f^{l_1}v_1\T\ldots\T f^{l_n}v_n,
\eean
where the sum is over all $\la$-admissible $\bar{l}\in\Zb^n_l$.

\begin{thm}\label{cite3'}
(Corollaries 17, 18 in \cite{MV1}.)
Let $p\in\C$ and $\la\in \C^n$ satisfy conditions \Ref{step}-\Ref{weights2} and  let $W\in\G^{sing}$ be a $\la$-admissible function. Then

1. The function $\Psi_W(z,\la)$ is a meromorphic solution of the qKZ equation with values in $(V_{\la_1}\T\ldots\T V_{\la_n})_l^{sing}$. 

2.  The function $\Psi^{\rm adm}_W(z,\la)$ is a meromorphic solution of the qKZ equation with values in $(L_{\la_1}\T\ldots\T L_{\la_n})_l^{sing}$. 
\end{thm}

The meromorphic solutions of the qKZ equation defined by \Ref{Psi} and \Ref{Psi adm} are called the \emph{hypergeometric
solutions}.

\subsection{Relations with representation theory}
Let $p\in\C$, $z,\la\in\C^n$ satisfy conditions \Ref{step}-\Ref{weights3}.

The weight space $(V_{\la_1}\T\ldots\T V_{\la_n})_l^*$ is identified with the hypergeometric rational space $\F(z,\la)$ by the map 
\be
{\frak a}(z,\la):\; (f^{l_1}v_1\T\ldots\T f^{l_n}v_n)^* \mapsto w_{\bar{l}}(t,z,\la).
\ee
Here $\{(f^{l_1}v_1\T\ldots\T f^{l_n}v_n)^*\,|\,l_1,\dots,l_n\in\Z_{\ge 0}\}$ is the basis of $(V_{\la_1}\T\ldots\T V_{\la_n})^*$, dual to the standard basis of $V_{\la_1}\T\ldots\T V_{\la_n}$ given by $\{f^{l_1}v_1\T\ldots\T f^{l_n}v_n\,|\,l_1,\dots,l_n\in\Z_{\ge 0}\}$, see Lemma 4.5, Corallary 4.8 in \cite{TV1}.

For  $\bar{l}\in\Zb^{n}_l$, define the \emph{weight coefficient} $c_{\bar{l}}(\la)$ by
\be
c_{\bar{l}}(\la)=
\prod_{m=1}^n\prod_{s=0}^{l_m-1}
\frac{\sin(\pi(s+1)/p)\sin(\pi(2\la_m-s)/p)}{\sin(\pi/p)}.
\ee

Let $q=e^{\pi i/p}$. The weight  space $(V_{\la_1}^q\T\ldots\T V_{\la_n}^q)_l$ is identified with the trigonometric hypergeometric space $\G (z,\la)$ by the map
\bean\label{identify}
{\frak b}(z,\la) :\; f_q^{l_1}v_1^q\T\ldots\T
f_q^{l_n}v_n^q \mapsto 
c_{\bar{l}}(\la)W_{\bar{l}}(t,z,\la),
\eean
see Lemma 4.17 and Corallary 4.20 in \cite{TV1}. The subspace of singular vectors $(V_{\la_1}^q\T\ldots\T V_{\la_n}^q)_l^{sing}$ is identified with the singular hypergeometric space,
\be
{\frak b}(z,\la)((V_{\la_1}^q\T\ldots\T V_{\la_n}^q)_l^{sing})=\G^{sing}(z,\la),
\ee
see Corallary 4.21 in \cite{TV1}. Moreover, the map ${\frak b}(z,\la)$ is factorized to an isomorphism 
\be
{\frak b}(z,\la)\,:\, (L_{\la_1}^q\T\ldots\T L_{\la_n}^q)_l\to \G_{\rm adm}(z,\la),
\ee
given by the same formula \Ref{identify}. The subspace of singular vectors is identified with the $\la$-admissible singular hypergeometric trigonometric space,
\be
{\frak b}(z,\la)((L_{\la_1}^q\T\ldots\T L_{\la_n}^q)_l^{sing})=\G_{\rm adm}^{sing}(z,\la).
\ee

We have the hypergeometric map,
\be
s(z,\la):\;
(V^q_{\la_1}\T\ldots\T V^q_{\la_n})_l^{sing}
\rightarrow
(V_{\la_1}\T\ldots\T V_{\la_n})_l^{sing},
\ee
defined for $v^q\in (V^q_{\la_1}\T\ldots\T V^q_{\la_n})_l^{sing}$ by 
$v^q\,\mapsto \,\Psi_{{\frak b}(z,\la)v^q}(z,\la)$, cf. Theorem 14 in \cite{MV1}, \cite{TV1}.

This map is factorized to the hypergeometric map,
\be
s(z,\la):\;
(L^q_{\la_1}\T\ldots\T L^q_{\la_n})_l^{sing}
\rightarrow
(L_{\la_1}\T\ldots\T L_{\la_n})_l^{sing},
\ee
defined for $v^q\in (L^q_{\la_1}\T\ldots\T L^q_{\la_n})_l^{sing}$ by 
$v^q\,\mapsto \,\Psi^{\rm adm}_{{\frak b}(z,\la)v^q}(z,\la)$,
cf. \cite{MV1}.

\begin{thm} (Theorem 24 in \cite{MV1}.) Let $p,z,\la$ satisfy conditions \Ref{step}-\Ref{weights3}. Let also $\Sum_{m=1}^n2\la_m-2l+k+1\not\in p\Z_{<0}$ for $k=1,\dots,l-1$. Then the map
$s(z,\la)\,:\,(L^q_{\la_1}\T\ldots\T L^q_{\la_n})_l^{sing}\rightarrow (L_{\la_1}\T\ldots\T L_{\la_n})_l^{sing}$ is an isomorphism of vector spaces.
\end{thm}

In this paper we describe the image and kernel of the hypergeometric map $s(z,\la)$ when the resonance condition, $\Sum_{m=1}^n2\La_m-2l+p+k+1=0$, holds for some $k\in\{1,\dots,l-1\}$.

\section{The hypergeometric solutions in the case of resonance}\label{s(z)}
\subsection{The space of quantized conformal blocks}
Let $L_{\la_i}$ be the irreducible $sl(2)$ module with highest weight $\la_i\in\C$, $i=1,\dots,n$. Let
\bean\label{e(z)-Yangian}
e(z)=T_{21}^{(2)}-T_{11}^{(1)}T_{21}^{(1)}\in Y(gl(2)).
\eean
Then, for each $\la$-admissible $\bar{l}\in\Z_{\ge 0}^n$, we have
\bean\label{e(z)}
e(z)f^{l_1}v_1\T\ldots\T f^{l_n}v_n=\sum_{j=1}^n \left(z_j+h^{(j)}+\sum_{s=j+1}^n2h^{(s)}\right)e^{(j)}f^{l_1}v_1\T\ldots\T f^{l_n}v_n=
\eean
\be
\sum_{j=1}^n (2\la_i-m_i+1)l_i(z_j+\la_j-l_j+\sum_{s=j+1}^n2(\la_s-l_s))f^{l_1}v_1\T\ldots\T f^{l_j-1}v_j\T\ldots\T f^{l_n}v_n,\notag
\ee
where $h^{(s)}, e^{(s)}$ are the elements $h,e\in sl(2)$ acting in the $s$-th factor.

Define the \emph{space of quantized conformal blocks $C(z)$}. Let
\bean\label{conformal blocks}
C(z)=\{m\in (L_{\la_1}\T\ldots\T L_{\la_n})^{sing}_l\,|\, (e(z))^km=0\},
\eean
if the resonance condition,
\bean\label{resonance condition} 
\Sum_{i=1}^n2\la_i-2l+p+k+1=0,
\eean
holds for some $k\in\N$, and $C(z)=(L_{\la_1}\T\ldots\T L_{\la_n})^{sing}_l$ otherwise.

Note that if $k>l$ then $C(z)=(L_{\la_1}\T\ldots\T L_{\la_n})^{sing}_l$.

The subspace $(L_{\la_1}\T\ldots\T L_{\la_n})^{sing}_l$ is called a \emph{resonance subspace} if the resonance condition \Ref{resonance condition}  holds for some $k\in\N$, $k\le l$.

For $x\in\C$, let
\bean\label{e(z)+p}
e_x(z)=e(z)+xT_{21}^{(1)}\in Y(gl(2)).
\eean
If \Ref{resonance condition} holds for some $k\in\N$, then the space $\{m\in (L_{\la_1}\T\ldots\T L_{\la_n})^{sing}_l\,|\, (e_x(z))^km=0\}$ does not depend on $x$ and coincides with the space of conformal blocks $C(z)$, see \cite{MV2}.

\begin{thm}\label{subbundle} (Theorem 2 in \cite{MV2}.)
The space of conformal blocks $C(z)$ is invariant with respect to the qKZ connection,  
\be
K_i(z)C(z)=C(z_1,\dots,z_i+p,\dots,z_n),
\ee
as well as with respect to permutations of variables, 
\be
PR_{M_iM_{i+1}}(z_i-z_{i+1})C(z)=C(z_1,\dots,z_{i+1},z_i,\dots,z_n).
\ee
\end{thm}

\subsection{The image of $s(z,\la)$.}
\begin{thm}\label{image1} Let $p,z,\la$ satisfy conditions \Ref{step}-\Ref{weights3} and let $L_{\la_i}$be the irreducible $sl(2)$ module with highest weight $\la_i$, $i=1,\dots,n$. Then the image of the map $s(z,\la)$ belongs to the space of conformal blocks $C(z)$, i.e. every hypergeometric solution \Ref{Psi adm} takes values in the space of conformal blocks.
\end{thm}

{\bf Example.} Let $l=k=1$, $n=3$. Then the solutions of the qKZ equation have the form 
\be
\Psi^{\pm}(z,\la)=I_1^{\pm}(z,\la)\,fv_1\T v_2\T v_3 +I_2^{\pm}(z,\la)\,v_1\T fv_2\T v_3+ I_3^{\pm}(z,\la)\,v_1\T v_2\T fv_3,
\ee
where
\be
I_m^\pm(z,\la)=\int \Phi(t,z,\la)w_m(t,z,\la)\left(\prod_{j=1}^3 \sin(\pi(t-z_j-\la_j)/p)\right)^{-1} e^{\pm \pi it/p}dt,
\ee 
$m=1,2,3$, and the weight functions $w_m(t,z,\la)$ are given by \Ref{one point}. According to Theorem~\ref{image1}, if $p+2\sum\limits_{j=1}^3\la_j=0$ then the coordinate functions $I_m^\pm(z,\la)$ satisfy the algebraic equation,
\be
\sum_{i=1}^32\la_i\left(z_i+\la_i+2\sum_{j=i+1}^3\la_j\right)I^\pm_i(z,\la)=0.
\ee

\begin{proof}
Theorem~\ref{image1} is trivial for non-resonance subspaces $(L_{\la_1}\T\ldots\T L_{\la_n})^{sing}_l$. Let $(L_{\la_1}\T\ldots\T L_{\la_n})_l^{sing}$ be a resonance subspace, so that $\Sum_{i=1}^n2\la_i-2l+p+k+1=0$ holds for some $k\in\N$, $k\le l$. Let $W\in\G^{sing}$ be an admissible trigonometric weight function. We have to prove that $\Psi^{\rm adm}_W(z,\la)\in C(z)$, where $\Psi^{\rm adm}_W(z,\la)$ is given by formula \Ref{Psi}. Let $e_p(z)\in Y(gl(2))$ be given by \Ref{e(z)+p}. We have to prove that $(e_p(z))^k\Psi^{\rm adm}_W(z,\la)=0$, i.e all coordinates of this vector equal zero. It is equivalent to the statement that for each $\la$- admissible $\bar{m}\in\Zb^n_{l-k}$, we have
\be
I\left({\frak a}(z,\la)(f_p(z))^k (f^{m_1}v_1\T\ldots\T f^{m_n}v_n)^*, \,W\right)=0,
\ee
where $f_p(z)\in\End((L_{\la_1}\T\ldots\T L_{\la_n})^*)$ is defined for $\la$-admissible $\bar{m}\in Z_{\ge 0}^n$, by
\be
f_p(z)(f^{m_1}v_1\T\ldots\T f^{m_n}v_n)^*=
\left.\sum_{i=1}^n\right((m_i+1)(2\la_i-m_i)\times
\ee
\be
\left.
\times(z_i+p+\la_i-m_i+\sum_{j=i+1}^n2(\la_j-m_j))\right)(f^{m_1}v_1\T\ldots\T f^{m_i+1}v_i\T\ldots\T f^{m_n}v_n)^*,
\ee
cf. formula \Ref{e(z)}.
 
For a function $\phi(t)=\phi(t_1,\dots,t_l)$, define the \emph{discrete partial derivatives} $(D_i\phi)(t)$, $i=1,\dots,l$, by
\be
(D_i\phi)(t)=\phi(t_1,\dots,t_i+p,\dots,t_l)-\phi(t).
\ee

Recall that $w_{\bar{l}}=\Sum_{\sigma\in{\Bbb S}^l}[\eta_{\bar{l}}]_\sigma^{rat}$, where $\eta_{\bar{l}}$ is defined in \Ref{eta}.

For fixed $z,\la\in\C^n$ and $p\in\C$, introduce a linear operator $\tilde f_p(z)$ acting in the space spanned over $\C$ by all functions $\eta_{\bar{l}}(t,z,\la)$, $\bar{l}\in\Z^n_{\ge 0}$, by the formula 
\be
\tilde{f}_p(z)\eta_{\bar{m}}=
\sum_{i=1}^n(m_i+1)(2\la_i-m_i)(z_i+p+\la_i-m_i+\sum_{j=i+1}^n2(\la_j-m_j))\eta_{(m_1,\dots,m_i+1,\dots,m_n)}.
\ee

\begin{lemma}\label{derivative} Let $\bar{m}\in\Zb^n_{l-k}$, then
\be
\sum_{\sigma\in{\Bbb S}}^l\left[\sum_{i=0}^k\left({k \choose i}\left(\prod_{j=1}^{k-i}(p+\sum_{i=1}^n2(\la_i-l)+k+j)\right)\right)((\tilde{f}_p(z))^i\eta_{\bar{m}}) (t_{k+1-i},\dots,t_l)\right]_\sigma^{rat}=
\ee
\bean\label{derivative formula}
=\sum_{\sigma\in{\Bbb S}}^l[(\Phi(t,z,\la))^{-1}\mu_{\bar{m}}(t,z,\la)]_\sigma^{rat},
\eean
where 
\be
%\mu_{\bar{m}}(t,z,\la)=(D_1(\phi_1(t,z,\la)(D_2(\phi_2(t,z,\la)\ldots(D_k(\phi_k(t,z,\la)
\mu_{\bar{m}}(t,z,\la)=(D_1(\phi_1(D_2(\phi_2\ldots(D_k(\phi_k
\eta_{\bar{m}}(t_{k+1},\dots,t_l,z,\la)))\dots)))),
\ee
\be
\phi_a(t,z,\la)=
t_a\prod_{i=1}^n
\frac{\Gamma((t_a-z_i+\la_i)/p)}{\Gamma((t_a-z_i-\la_i)/p)}
\prod_{j=a+1}^l
\frac{\Gamma((t_a-t_j-1)/p)}{\Gamma((t_a-t_j+1)/p)}.
\ee
\end{lemma}

{\bf Example.} Let $l=k=1$. Then
\bea
\lefteqn{D_1 (t\Phi)=}
\\&&
\left(p+\sum_{i=1}^n 2\la_i\right)\Phi(t,z,\la)+\sum_{i=1}^n 2\la_i\left(z_i+p+\la_i+\sum_{j=i+1}^n 2\la_j\right)w_i(t,z,\la)\Phi(t,z,\la),
\eea
where $w_i(t,z,\la)$ are the weight functions given by \Ref{one point}.

\begin{proof}
Introduce the numbers $y_1,\dots,y_{n+l-1}$ by
\be
(y_1,\dots,y_{n+l-1})=(t_2,\dots,t_{l_1+1},z_1,t_{l_1+2},\dots,t_{l_1+l_2+1},z_2,t_{l_1+l_2+2},\dots,z_n),
\ee
and the numbers $\Dl_1,\dots,\Dl_{n+l-1}$ by
\be
(\Dl_1,\dots,\Dl_{n+l-1})=(-2,\dots,-2,2\la_1,-2,\dots,-2,2\la_2,-2,\dots,2\la_n).
\ee
Order the numbers $y_1,\dots,y_{n+l-1}$ by the rule: $y_i\prec y_j$ if and only if $i<j$.

For $i=1,\dots,n+l-1$, introduce
\be
\eta_i(t_1,y,\Dl)=\frac{1}{t_1-y_i-\Dl_i} \prod_{j=1}^{i-1}\frac{t_1-y_j+\Dl_j}{t_1-y_j-\Dl_j},
\ee
cf. formula \Ref{one point}.

For $i=1,\dots, n+l-1$, we have the following equalities, proved by induction on $i$,
\be
\prod_{s=1}^i\frac{t_1-y_s+\Dl_s}{t_1-y_s-\Dl_s}-1=\sum_{s=1}^i 2\Dl_s\eta_s(t_1,y,\Dl),
\ee
\be
t_1\eta_i(t_1,y,\Dl)-1=\sum_{j=1}^{i-1}2\Dl_j\eta_j(t_1,y,\Dl)+(y_i+\Dl_i)\eta_i(t_1,y,\Dl).
\ee

For $\bar{l}\in\Zb^n_{l-1}$, $j=1,\dots,n+l-1$, we have
\be
\sum_{\sigma\in{\Bbb S}^l}[\eta_{\bar{l}}(t_2,\dots,t_l,z,\la)\eta_j(t_1,y,\Dl)]_\sigma^{rat}=
\sum_{\sigma\in{\Bbb S}^l}[a_j\eta_{(l_1,\dots,l_i+1,\dots,l_n)}(t,z,\la)]_\sigma^{rat},
\ee
where $i\in\{1,\dots,n\}$ is such that $z_{i-1}\prec y_j\preceq z_{i}$; $a_j=l_i+1$ if $y_j=z_j$ and
$a_j=\half(l_i+1)$ otherwise.

Let $j\in\{1,\dots,n+l-1\}$ be such that $z_{i-1}\prec y_j=t_s\prec z_i$. Then for $\bar{l}\in\Zb^n_{l-1}$, $j=1,\dots,n+l-1$, we have
\be
\sum_{\sigma\in{\Bbb S}^l}[\eta_{\bar{l}}(t_2,\dots,t_l,z,\la)y_j\eta_j(t_1,y,\Dl)]_\sigma^{rat}=
\sum_{\sigma\in{\Bbb S}^l}[(z_i+\la_i)\eta_{\bar{l}}(t,z,\la)\eta_j(t_1,y,\Dl)]_\sigma^{rat}.
\ee

Using the above formulas, we compute
\be
\sum_{\sigma\in{\Bbb S}^l}\left[(t_1+p)\eta_{\bar{l}}(t_2,\dots,t_l,z,\la)\left(\prod_{s=1}^{n+l-1}\frac{t_1-y_s+\Dl_s}{t_1-y_s-\Dl_s}-1\right)\right]_\sigma^{rat}=\sum_{i=1}^n2\La_i-2(l-1)+
\ee
\be
+\sum_{i=1}^n(l_i+1)(2\la_i-l_i)\left(z_i+\la_i-l_i+\sum_{j=i+1}^n2(\la_i-l_i)\right)\sum_{\sigma\in{\Bbb S}^l}[\eta_{(l_1,\dots,l_i+1,\dots,l_n)}(t,z,\la)]_\sigma^{rat}.
\ee
This proves the Lemma when $k=1$.
Lemma~\ref{derivative} for other $k$ is proved similarly by induction on k.
\end{proof}

We deduce the Theorem from Lemma~\ref{derivative}. Under resonance condition \Ref{resonance condition}, the terms in the left hand side of \Ref{derivative formula} corresponding to $i=1,\dots,k-1$ equal zero, and the only non-zero term (corresponding to $i=k$) equals ${\frak a}(z,\la)(f(z))^k (f^{m_1}v_1\T\ldots\T f^{m_n}v_n)^*$. So, 
\be
I\left({\frak a}(z,\la)(f(z))^k (f^{m_1}v_1\T\ldots\T f^{m_n}v_n)^*,W\right)(z,\la)=
\ee
\be
\int \Phi(t,z,\la)\left( \sum_{\sigma\in{\Bbb S}}^l[(\Phi(t,z,\la))^{-1}\mu_{\bar{m}}(t,z,\la)]_\sigma^{rat}\right) W d^lt=
\ee
\be
l! \int \Phi(t,z,\la)(\Phi(t,z,\la))^{-1}\mu_{\bar{m}}(t,z,\la)) W d^lt=l! \int \mu_{\bar{m}}(t,z,\la)) W d^lt,
\ee
where the integration is over a suitable contour, see \cite{MV1}, and in the second equality we use the symmetry of the initial integrand, see the Remark in Section~\ref{actions}.

The function $W(t,z,\la)$ is $p$-periodic with respect to each of $t_1,\dots,t_l$. According to Lemma~\ref{derivative}, $\mu_{\bar{m}}(t,z,\la)=D_1g_{\bar{m}}(t,z,\la)$ for some function $g_{\bar{m}}(t,z,\la)$. Hence $\mu_{\bar{m}}(t,z,\la) W(t,z,\la)=D_1(g_{\bar{m}}W)(t,z,\la)$. Note that the function $g_{\bar{m}}(t,z,\la)$ has a polynomial growth as $t_j$ goes to $\pm i\infty$. All functions $W(t,z,\la)\in\G^{sing}(z,\la)$ decay exponentially as $t_j$ goes to $\pm i\infty$. Therefore, all the integrals converge.

If $\la_i\ll 0$ for all $i=1,\dots,n$, then the integral $\int\limits_{\Real t_i=0,\atop i=1,\dots,l} D_1(g_{\bar{m}}W)(t,z,\la)d^lt$ is zero, since the function $g_{\bar{m}}(t,z,\la)W(t,z,\la)$ has no poles between $\{t\in\C^l\,|\,\Real t_i=0,\, i=1,\dots,l\}$ and $\{t\in\C^l\,|\, \Real t_1=\Real p,\,\Real t_i=0,\, i=2,\dots,l\}$. The integral is a meromorphic function of $\la$, hence it is zero for all $\la$, cf. Lemma 9.5 in \cite{TV1}.
\end{proof}

\subsection{The kernel of $s(z,\la)$.}

Let $p,z,\la$ satisfy conditions \Ref{step}-\Ref{weights3}.
Let $L_{\la_i}$ be the irreducible $sl(2)$ module with highest weight $\la_i\in\C$ and let $L^q_{\la_i}$ be the corresponding $U_qsl(2)$ module with highest weight $q^{\la_i}$, $i=1,\dots,n$.

Let $(L_{\la_1}\T\ldots\T L_{\la_n})_l^{sing}$ be a resonance subspace, so that $\Sum_{i=1}^n2\la_i-2l+p+k+1=0$ holds for some $k\in\N$, $k\le l$.
\begin{lemma}\label{f**k}
Let $v^q\in (L_{\la_1}^q\T\ldots\T L_{\la_n}^q)_{l-k}^{sing}$. Then $(f_q)^kv^q\in(L_{\la_1}^q\T\ldots\T L_{\la_n}^q)_l^{sing}$.
\end{lemma}
\begin{proof}
We have
\be
e_q(f_q)^kv^q=[\sum_{i=1}^n2\la_i-2(l-k)-k+1]_q[k]_q(f_q)^{k-1}v^q,
\ee
where $[k]_q=(q^k-q^{-k})/(q-q^{-1})$. Recall that $q=e^{\pi i /p}$, so $[\Sum_{i=1}^n2\la_i-2(l-k)-k+1]_q=[-p]_q=0$.
\end{proof}

\begin{thm}\label{kernel1}
Let $v^q\in(L^q_{\la_1}\T\ldots\T L_{\la_n}^q)_l^{sing}$ have the form $v^q=(f_q)^k\tilde{v}^q$ for some $\tilde{v}^q\in(L^q_{\la_1}\T\ldots\T L^q_{\la_n})_{l-k}^{sing}$. Then the hypergeometric solution $s(z,\la)\,v^q$ equals zero.
\end{thm}
\begin{proof}
Let $W={\frak b}(z,\la)v^q$. We need to prove that $\Psi^{\rm adm}_W(z)$ given by \Ref{Psi adm} is zero.

The space $\bigoplus\limits_{l=0}^\infty \G_l^n(z,\la)$ has a $U_qsl_2$ module structure such that the map ${\frak b}(z,\la)$ is an intertwiner of $U_qsl(2)$ modules, see \cite{TV1}. The action of $U_qsl(2)$ in the space $\bigoplus\limits_{l=0}^\infty \G_l^n(z,\la)$ is given by formulas (4.16) in \cite{TV1}. In particular, for any $X(t,z,\la)\in\G_l^n(z,\la)$, we have
\be
(f_qX)(t_1,\dots,t_{l+1})=
\left.\left.\exp\left(-\pi i(l+\sum_{m=1}^n\la_m)/p\right)\sum_{a=1}^{l+1}\right[X(t_2,\dots,t_{l+1})
\right(\exp(2\pi il/p)\times
\ee
\be
\left.\left.
\times\prod_{m=1}^n\frac{\sin(\pi(t_1-z_m+\la_m)/p)}{\sin(\pi(t_1-z_m-\la_m)/p)}
\prod_{b=2}^{l+1}\frac{\sin(\pi(t_1-t_b-1)/p)}{\sin(\pi(t_1-t_b+1)/p)}-
\exp(2\pi i\sum_{m=1}^n\la_m/p)\right)\right]_{(1,a)}^{trig},
\ee
where $(1,a)\in {\Bbb S}^{l+1}$ are transpositions.

Therefore, $W={\frak b}(z,\la)(f_q)^k\tilde{v^q}$ has the form $W=\Sum_{\sigma\in{\Bbb S}^l}[\Upsilon]_\sigma^{trig}$, where 
\bean\label{Upsilon}
\Upsilon(t_1,\dots,t_l,z,\la)=\left.a(\la)X(t_{k+1},\dots,t_l)\prod_{i=1}^k\right(\exp\left(2\pi i(l-i)/p\right)\times
\eean
\be
\left.
\times\prod_{m=1}^n\frac{\sin(\pi(t_i-z_m+\la_m)/p)}{\sin(\pi(t_i-z_m-\la_m)/p)}
\prod_{b=i+1}^l\frac{\sin(\pi(t_i-t_b-1)/p)}{\sin(\pi(t_i-t_b+1)/p)}-
\exp(2\pi i\sum_{m=1}^n\la_m/p)\right),
\ee
for some $X\in\G^{sing}_{l-k}(z,\la)$ and some holomorphic function $a(\la)$.

Fix $w=w_{\bar{l}}$, where $\bar{l}\in\Zb^n_l$. 
Fix $\la$ such that $\Imag \la_i\neq 0$, and $\Sum_{i=1}^n2\la_i-2l+p+k+1=0$. We need to prove that $I(w,W)(z,\la)=0$, where $I(w,W)$ is the hypergeometric integral defined in Section~\ref{integrals}. The integrand $\Phi wW$ is a meromorphic function of $t,z,\la$ with simple poles located at most at
\be
t_i=z_j\pm(\la_j+kp),\qquad k\in\Z_{\ge0},\qquad i=1,\dots,l,\qquad j=1,\dots,n,
\ee
and at 
\bean\label{poles3}
t_i-t_j=\pm(1-kp), \qquad k\in\Z_{\ge 0},\qquad i,j=1,\dots,l,\qquad i<j.
\eean

Let $z\in\C^n$ be such that $|\Imag\la_i|\ll\Imag z_1\ll\ldots\ll\Imag z_n$. Let $p$ be a real negative number. Let $D(z,\la)\subset\C$ be a curve such that $D(z,\la)$ coincides with the imaginary line at a neighbourhood of infinity, $D(z,\la)$ consists of straight non-horisontal segments, the points 
$z_j-\la_j$, $j=1,\dots,n$, are to the right of $D(z,\la)$ and the points $z_j+\la_j$, $j=1,\dots,n$, are to the left of $D(z,\la)$, see the picture.

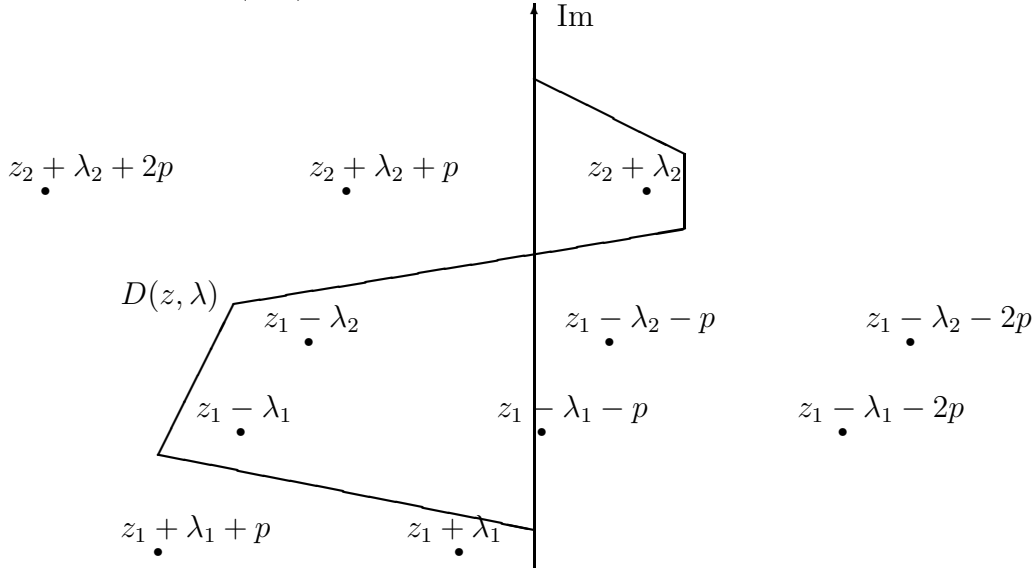
\setlength{\unitlength}{1cm}
\begin{picture}(14,7.5)
\thinlines
\put (7,0){\vector(0,1){7.5}}
\put (7.3,7.2){Im}
\put (1.5,3.5){$D(z,\lambda)$}
\thicklines
\put (7,6.5){\line(0,1){0.7}}
\put (7,6.5){\line(2,-1){2}}
\put (9,4.5){\line(0,1){1}}
\put (3,3.5){\line(6,1){6}}
\put (2,1.5){\line(1,2){1}}
\put (2,1.5){\line(5,-1){5}}
\put (7,0){\line(0,1){0.5}}

\put (4,3){\circle*{0.1}}
\put (3.4,3.2){$z_1-\lambda_2$}
\put (8,3){\circle*{0.1}}
\put (7.4,3.2){$z_1-\lambda_2-p$}
\put (12,3){\circle*{0.1}}
\put (11.4,3.2){$z_1-\lambda_2-2p$}

\put (3.1,1.8){\circle*{0.1}}
\put (2.5,2){$z_1-\lambda_1$}
\put (7.1,1.8){\circle*{0.1}}
\put (6.5,2){$z_1-\lambda_1-p$}
\put (11.1,1.8){\circle*{0.1}}
\put (10.5,2){$z_1-\lambda_1-2p$} 

\put (8.5,5){\circle*{0.1}}
\put (7.7,5.2){$z_2+\lambda_2$}
\put (4.5,5){\circle*{0.1}}
\put (4,5.2){$z_2+\lambda_2+p$}
\put (0.5,5){\circle*{0.1}}
\put (0.0,5.2){$z_2+\lambda_2+2p$}

\put (6,0.2){\circle*{0.1}}
\put (5.3,0.4){$z_1+\lambda_1$}
\put (2,0.2){\circle*{0.1}}
\put (1.5,0.4){$z_1+\lambda_1+p$}

\end{picture}
\bigskip

The analytic continuation of the hypergeometric integral $I(w,W)$ with respect to $\la$ gives
\be
I(w,W)=\int\limits_{t_i\in D(z,\la)\atop i=1,\dots,l}\Phi wWd^lt,
\ee
cf. \cite{MV1}.

Due to the symmetry of the function $\Phi wW$, it is enough to prove that the integral of $\Phi w\Upsilon$ equals zero, $I(w,\Upsilon)=0$, see the Remark in Section~\ref{actions}.

Let $A,B$ be large positive real numbers. Let $D_A(z,\la)$ be the curve $D(z,\la)$ truncated at infinity, namely,
$D_A(z,\la)=\{x\in D(z,\la)\,|\, |\Imag x|\le A\}$. Let
\be
I_{AB}(w,\Upsilon)=\int\limits_{t_j\in D_B(z,\la)\atop j=k+1,\dots,l}\int\limits_{t_i\in D_A(z,\la)\atop i=1,\dots,k}\Phi w\Upsilon d^lt.
\ee
We have $I_{AB}(w,\Upsilon)\to I(w,\Upsilon)$ when $A,B\to +\infty$.

Let $D_A^+=\{x\in\C\,|\, \Real x\ge 0,\,|x|=A\}$ and $D_A^-=\{x\in\C\,|\, \Real x\le 0\,|x|=A\}$ are halves of the circle of radius $A$ centered at the origin.  

The function $\Upsilon(t,z,\la)$ is given by \Ref{Upsilon} as a product of $k$ factors, and each of the factors is a sum of two terms. Represent the function $\Upsilon(t,z,\la)$ as a sum of $2^k$ terms, $\Upsilon=\Sum_{s=1}^{2^k}\Upsilon_{(s)}$. 
Let $E_{(s)}\subset\{1,\dots,k\}$ consist of the indices $i\in\{1,\dots,k\}$ such that $\Upsilon_{(s)}$ contains the factor
\be
\prod_{m=1}^n\frac{\sin(\pi(t_i-z_m+\la_m)/p)}{\sin(\pi(t_i-z_m-\la_m)/p)}.
\ee
The function $\Phi w\Upsilon_{(s)}$ has no pole at
\be
t_i=z_j-(\la_j+kp),\qquad k\in\Z_{\ge0},\qquad i\in E_{(s)},\qquad j=1,\dots,n,
\ee
\be
t_i=z_j+(\la_j+kp),\qquad k\in\Z_{\ge0},\qquad i\in\{1,\dots,k\},\,i\not\in E_{(s)}, \qquad j=1,\dots,n.
\ee
Consider the integral $I_{AB}(w,\Upsilon_{(s)})$. We move the contours of integration with respect to $t_i$, $i=1,\dots,k$, to the right if $i\not\in E_{(s)}$ and to the left if $i\in E_{(s)}$. It easy to see that we do not encounter poles \Ref{poles3}, therefore, we get
\be
I_{AB}(w,\Upsilon_{(s)})=\int\limits_{t_a\in D_B(z,\la)\atop a=k+1,\dots,l}\int\limits_{t_i\in D_A(z,\la)\atop i=1,\dots,k}\Phi w\Upsilon_{(s)}d^lt=
\int\limits_{t_a\in D_B(z,\la)\atop a=k+1,\dots,l}\int\limits_{t_j\in D^+_A(z,\la)\atop j\not\in E_{(s)}}\int\limits_{t_i\in D^-_A(z,\la)\atop i\in E_{(s)}}\Phi w\Upsilon_{(s)}d^lt.
\ee

Now, we fix $t_{k+1},\dots,t_l$ and estimate the integrand using the Stirling formula as $A\to+\infty$. We have
\be
\Phi (t,z,\la)\Upsilon_{(s)}(t,z,\la)\,\sim\, c_1A^{k\Sum_{i=1}^n2\la_i/p-2(l-k)k/p-k(k-1)/p},\qquad w(t,z,\la)\,\sim\,c_2A^{-k},
\ee
where $c_1,c_2\in\C$ are some quantities independent on $A$. Multiplying and using the resonance condition, we get
\be
\Phi w\Upsilon_{(s)}\,\sim\,c_1c_2A^{k/p(\sum_{i=1}^n2\la_i-2l+k+1)-k}=c_1c_2A^{-2k}.
\ee
We integrate with respect to each of $t_1,\dots,t_k$ over a semi-circle of radius $A$, hence the integral with respect to $t_1,\dots,t_k$ can be estimated by
\be
\int\limits_{t_j\in D^+_A(z,\la)\atop j\not\in E_{(s)}}\int\limits_{t_i\in D^-_A(z,\la)\atop i\in E_{(s)}}\Phi w\Upsilon_{(s)}d^lt\,\sim\, c_1c_2A^{-2k}(\pi A)^k=c_3A^{-k},
\ee
for some $c_3\in\C$ independent on $A$. Therefore, $I_{AB}(w,\Upsilon)\to 0$ as $A\to +\infty$. Thus, $I(w,\Upsilon)(z,\la)=I(w,W)(z,\la)=0$ for our choice of $z,\la$. 

The integral $I(w,W)$ is a meromorphic function of $z,\la$, therefore, $I(w,W)(z,\la)=0$ for all $z,\la\in\C^n$.

This proof can be easily adjusted to any $p\in\C$ with properties \Ref{step},\Ref{step3}.
\end{proof}

\subsection{Finite dimensional representations and an integer step}

Let $\la_i\in\La^+$, $i=1,\dots,n$, are non-negative half-integers. Let $-p,l\in\Z_{>0}$, and
\bean\label{good condition}
2\la_i\le -p-2,\qquad i=1,\dots,n.
\eean
Note that conditions \Ref{step}-\Ref{weights2} and \Ref{step3}-\Ref{weights3} are automatically satisfied.

Let
\bean\label{bad condition}
l\le 2\la_i,\qquad i=1,\dots,n.
\eean

\begin{thm}\label{kernel2} Let $p,z,\la$ satisfy conditions \Ref{step}-\Ref{weights3}, \Ref{good condition} and \Ref{bad condition}. Let the resonance condition \Ref{resonance condition} hold for some $k\in\N$. Let $L_{\la_i}$ be the finite dimesional irreducible module with highest weight $\la_i$ and let $L_{\la_i}^q$ be the corresponding $U_qsl(2)$ module with highest weight $q^{\la_i}$, $i=1,\dots,n$. Then for $v^q\in(L_{\la_1}^q\T\ldots\T L^q_{\la_n})_l^{sing}$, the hypergeometric solution $s(z,\la)v^q$ equals zero if and only if $v^q$ has the form $v^q=(f_q)^k\tilde{v}^q$ for some $\tilde{v}^q\in(L^q_{\la_1}\T\ldots\T L^q_{\la_n})_{l-k}^{sing}$. 
\end{thm}

\begin{proof}
 The determinant of the hypergeometric pairing is given by formula (5.15) in \cite{TV1}. It follows from formula (5.15) in \cite{TV1} that under conditions of the Theorem the dimension of the kernel of the map $s(z,\la)$ is not greater than
\be
{n+l-k-2 \choose n-2}=\dim(V^q_{\la_1}\T\ldots\T V^q_{\la_n})_{l-k}^{sing}=
\dim(L_{\la_1}^q\T\ldots\T L^q_{\la_n})_{l-k}^{sing}.
\ee
Under conditions of Theorem~\ref{kernel2} the operator $(f_q)^k$ acting on $(L^q_{\la_1}\T\ldots\T L^q_{\la_n})_{l-k}^{sing}$ is non-degenerate, so, Theorem~\ref{kernel2} follows from Theorem~\ref{kernel1} and Lemma~\ref{f**k}.
\end{proof}

\begin{corollary}\label{corol} Under assumptions of Theorem~\ref{kernel2}, the space of hypergeometric solutions is naturally identified with the space
\be
(L_{\la_1}^q\T\ldots\T L^q_{\la_n})_l^{sing}/(f_q)^k((L_{\la_1}^q\T\ldots\T L^q_{\la_n})_{l-k}^{sing}).
\ee
\end{corollary}

Corollary~\ref{corol} follows from Theorem~\ref{kernel2}.

\begin{thm}\label{image2} 
Let $p,z,\la$ satisfy conditions \Ref{step}-\Ref{weights3}, \Ref{good condition} and \Ref{bad condition}. Let the resonance condition \Ref{resonance condition} hold for some $k\in\N$. Let $L_{\la_i}$ be the finite dimesional irreducible module with highest weight $\la_i$ and let $L_{\la_i}^q$ be the corresponding $U_qsl(2)$ module with highest weight $q^\la$, $i=1,\dots,n$. Then the image of the map $s(z,\la)\,:\,(L^q_{\la_1}\T\ldots\T L^q_{\la_n})_l^{sing}\to (L_{\la_1}\T\ldots\T L_{\la_n})_l^{sing}$  coincides with  the space of conformal blocks $C(z)$.
\end{thm}
\begin{proof}
Recall that for the differential KZ connection the space of conformal blocks $N(z)$ is defined as the kernel of the operator $(E(z))^k$ acting in a resonance subspace $(L_{\la_1}\T\ldots\T L_{\la_n})^{sing}_l$, where 
\be
E(z)f^{l_1}v_1\T\ldots f^{l_n}v_n=\sum_{i=1}^n l_i(2\la_i-l_i+1)z_i\,f^{l_1}v_1\T\ldots\T f^{l_i-1}v_i\T\ldots\T f^{l_n}v_n,
\ee
see [FSV1-3].
By Theorem 3.4 in \cite{Fn}, 
\be
\dim N(z)=\dim (L^q_{\la_1}\T\ldots\T L^q_{\la_n})_l^{sing}-\dim (L^q_{\la_1}\T\ldots\T L^q_{\la_n})_{l-k}^{sing}.
\ee

\begin{lemma}\label{clas>=quant} The dimensions  of the spaces of conformal blocks in the differential case and difference case are related by $\dim N(z)\ge\dim C(z)$.
\end{lemma}
\begin{proof}
For a function $\psi(z)$ with values in $(L_{\la_1}\T\ldots\T L_{\la_n})_{l-k}$, denote $\psi_{\bar{l}}(z)$ its $\bar{l}$-th coordinate function with respect to the basis $f^{l_1}v_1\T\ldots f^{l_n}v_n$, $\bar{l}\in\Zb^n_{l-k}$.

Let $r=\dim (L_{\la_1}\T\ldots\T L_{\la_n})^{sing}_l-\dim N(z)$ be the dimension of the image of the operator $(E(z))^k$ acting in $(L_{\la_1}\T\ldots\T L_{\la_n})^{sing}_l$. Then there exist vectors $w_1,\dots,w_r\in(L_{\la_1}\T\ldots\T L_{\la_n})^{sing}_l$ independent on $z$ such that for generic $z\in\C^n$, their images $(E(z))^kw_1,\dots,(E(z))^kw_r\in(L_{\la_1}\T\ldots\T L_{\la_n})_{l-k}$ are linearly independent. It means that there are indices $\bar{l}_1,\dots,\bar{l}_r\in\Zb^n_{l-k}$ such that the matrix $F(z)=\{(E(z))^kw_i)_{\bar{l}_j}\}_{j=1}^r$ is non-degenerate. Note that the determinant of $F(z)$ is a non-zero homogeneous polynomial of $z_1,\dots,z_n$ of degree $kr$. 

Consider the matrix $G=\{(e(z))^kw_i)_{\bar{l}_j}\}_{j=1}^r$. We have
\be
\det G(z)=\det F(z)+g(z),
\ee
where $g(z)$ is a polynomial in $z_1,\dots,z_n$ of degree less than $kr$. Therefore, $\det G(z)\neq 0$. 
\end{proof}

By Theorem~\ref{image1}, the image of the map $s(z,\la)$ belongs to the subspace of conformal blocks $C(z)$. By Theorem~\ref{kernel2}, the dimension of the image of the map $s(z,\la)$ is equal to the dimension of $N(z)$, which not smaller than the dimension of $C(z)$ by Lemma~\ref{clas>=quant}. Hence, $\dim C(z)=N(z)$ and the map $s(z,\la)$ is onto the subspace of conformal blocks $C(z)$.
\end{proof}

We expect that Theorems~\ref{kernel2} and \ref{image2} hold without the assumption $l\le 2\la_i$, $i=1,\dots,n$.

\begin{corollary}\label{final1} Under assumptions of Theorem~\ref{image2}, the dimensions  of the subspaces of conformal blocks in the differential case and difference case are the same $\dim N(z)=\dim C(z)$ and equal $\dim(L^q_{\la_1}\T\ldots\T L^q_{\la_n})^{sing}_l-\dim(L^q_{\la_1}\T\ldots\T L^q_{\la_n})^{sing}_{l-k}$.
\end{corollary}

Corollary~\ref{final1} follows from the proof of Theorem~\ref{image2}.

\begin{corollary}\label{final2}  Under assumptions of Theorem~\ref{image2}, the map $s(z,\la)$ is an isomorphism of vector spaces:
\be
s(z,\la)\,:\,(L^q_1\T\ldots\T L^q_n)^{sing}_l/f_q^k(L^q_1\T\ldots\T L^q_n)^{sing}_{l-k}\,\to\,C(z).
\ee
\end{corollary}

Corollary~{final2} follows from Theorems~\ref{kernel2} and \ref{image2}.

\end{document}